\documentclass[12pt,reqno]{amsart}

\usepackage{AKstyle3}
\usepackage[table]{xcolor}

\numberwithin{equation}{section}

\usepackage{caption}
\usepackage[labelfont=rm]{subcaption}

\usepackage[pagebackref=true, colorlinks]{hyperref}
\hypersetup{pdffitwindow=true,linkcolor=blue,citecolor=blue,urlcolor=blue,menucolor=blue}

\usepackage{comment}

\begin{document}

\author{Alex Kontorovich}
\thanks{Kontorovich is partially supported by
an NSF CAREER grant DMS-1455705, an NSF FRG grant DMS-1463940,  a BSF grant, a Simons Fellowship,  a von Neumann Fellowship at IAS, and the IAS's NSF grant DMS-1638352.}
\email{alex.kontorovich@rutgers.edu}
\address{Rutgers University, New Brunswick, NJ and Institute for Advanced Study, Princeton, NJ}

\author{Jeffrey   Lagarias}
\thanks{Lagarias is partially supported by NSF grants
 DMS-1401224 and DMS-1701576.}
\email{lagarias@umich.edu}
\address{University of Michigan, Ann Arbor, MI}

\title
{On Toric Orbits in the Affine Sieve}

\begin{abstract}
We give a detailed analysis of a heuristic model for the failure of ``saturation'' in instances of the Affine Sieve having toral Zariski closure. Based on this model, we formulate precise conjectures on several classical problems of arithmetic interest,  and test these against empirical data. 
\end{abstract}
\date{August 7, 2018}
\maketitle
\tableofcontents

\section{Introduction}\label{sec:intro}

The Fundamental Theorem of the Affine Sieve, introduced by Bourgain-Gamburd-Sarnak \cite{BourgainGamburdSarnak2010} and proved by Salehi Golsefidy-Sarnak \cite{SalehiSarnak2013} extends the Brun sieve to orbits of affine-linear group actions. 
The goal of this paper is to study the behavior of prime factors of orbits outside the 
 purview 
of this theorem.

More precisely, let $\G<\GL_N(\Q)$ be a finitely generated group,
that is, $\G = \langle A_1, A_2, \cdots, A_k \rangle$, 
 fix a base point $\bv_0\in\Q^N$, and let 
$$
\cO:=\G\cdot \bv_0 \ \subset \ \Z^N
$$
be the orbit of $\bv_0$ under $\G$, assumed to be integral\footnote{One can work more generally with entries in the ring of $S$-integers $\Z_S$, but we restrict to $\Z$ for ease of exposition. Note that there exist $\G<\GL_N(\Q)$
having no non-zero vector giving an integral orbit, e.g.,  $\G=\<A\>$ with $A=\mattwos\foh00{\foh}$.}.
Let $\gW(n)$ 
denote the number of primes dividing an integer $n$, counted with 
multiplicity.
Given $R\ge1$, 
an integer $n$ with $\gW(n)\le R$
is called 
{$R$-almost prime}. 
Fix a polynomial $f(x_1, x_2, ..., x_N) \in\Q[x_1,\dots,x_N]$ taking integer values on $\cO$, and let 
$$
\cO_R:=\{\bv\in\cO: \gW(f(\bv))\le R\}
$$
be the points in $\cO$ taking $R$-almost prime values under $f$. The pair $(\cO,f)$ is said to {\it saturate} if there exists some $R<\infty$ so that
\be\label{eq:sat}
\Zcl(\cO_R) = \Zcl(\cO).
\ee
Here $\Zcl$ refers to Zariski closure in affine space.\footnote{Recall that this Zariski closure can be thought of as the zero set of all polynomials vanishing on $\cO$.} 
The {\it saturation number} is the least $R$ for which \eqref{eq:sat} holds; this can be determined exactly or at least well-approximated in some special instances, see \cite{Kontorovich2014} for more discussion.
Let $V(f)$ be the 
affine
$\Q$-variety given by $f=0$.
In general, we assume that $f$ is non-constant on (any irreducible component of) $\Zcl(\cO)$. 
This is equivalent to 
\be\label{eq:Zcl}
\dim ( V(f) \cap \Zcl(\cO) ) < \dim  \Zcl(\cO),
\ee
 viewing the Zariski closure $\Zcl(\cO)$ inside $\C^{N}$. 
Then the aforementioned Fundamental Theorem
 of Salehi Golsefidy and Sarnak \cite[Theorem 1]{SalehiSarnak2013}, 
 states the following.
\begin{thm}[\cite{SalehiSarnak2013}]
Let $\G$ be a finitely generated subgroup of $GL_{N}(\Q)$ having Zariski closure
 $\bbG=\Zcl(\Gamma)$
in $GL_{N}(\C)$. Let ${\bv}_0 \in \Q^N$ and let  $\cO = \Gamma {\bv}_0 \subset \Z^N$ be the 
$\Gamma$-orbit of $\bv_0$.
Suppose  that $f(x) \in \Q[x_1,\dots,x_N]$ is such that $f(\cO)\subset\Z$ and \eqref{eq:Zcl} is satisfied.
Then the  pair $(\cO,f)$ saturates, as long as no algebraic torus\footnote{E.g. $(\C^\times)^n$.
} is a homomorphic image of the connected component $\bbG_0$ of the identity 
of  $\bbG$.
\end{thm}

In \cite[Appendix]{SalehiSarnak2013}, Salehi Golsefidy-Sarnak  give a heuristic argument, based on the Borel-Cantelli lemma, that the condition of having no tori is necessary in certain cases. 
Their model considered an algebraic torus  (that is, $\G$ is a free abelian group of rank $D$ with generators $A_1,\dots,A_D\in\GL_N(\Z)$, and there is a $g\in\GL_N(\C)$ so that for all $j$, the matrices $gA_j g^{-1}$ are diagonal)
and the test polynomial $f(x_{1,1}, ..., x_{N,N}) =\prod_{j=1}^k f_j(\bx)$, with $f_j(\bx)=j+\sum_{m,n=1}^N x_{m,n}^2 $. The test polynomial $f$  has (at least) $k$  irreducible
factors over $\Q[x_{1,1}, ... x_{N,N}]$, all of the same degree  (so they have roughly the same ``size'' on points of 
$\cO$). Their heuristic was
that the prime factorizations of the $k$ elements $f_j(\bx)$
evaluated at a point $\bx \in \cO$ ought to be ``independent," at least at the level of the number of
prime factors, $\gW(f_j (\bx))$, 
since they are just integer shifts of each other. 

In this paper, we refine this heuristic and make precise predictions on the failure of saturation in the toric case, which we then test empirically in a number of natural  settings of classical interest.

\subsection{Main Probabilistic Model}\

We model the $k$ irreducible factors of $f$ as $k$ randomly and independently chosen integers in an exponentially growing interval, depending on a parameter $n$.
  The parameter $n$ is to be viewed as modeling elements of a toral orbit, which grow exponentially.
  
\begin{thm}\label{thm:main1}
Let  $k\ge 1$ be a fixed integer. Fix a constant  $C >1$  and 
for each $n\ge1$, draw an integer  vector
$$
(x_{1,n},x_{2,n},\dots,x_{k,n}) \ \in \ [1,C^n]^k
$$
with uniform distribution. Then with probability one,
\be\label{eq:gWliminf}
\liminf_{n\ge1}{\gW(x_{1,n}\cdot x_{2,n}\cdots x_{k,n})\over \log n} \ = \ \gb_k,
\ee
where $\gb_k$ denotes the unique solution in $[0,k-1]$ to
\be\label{eq:gbkDef}
\gb_k(1-\log\gb_k+\log k) \ = \ k-1,
\ee
with $\beta_1=0$ and $\beta_k>0$ for $k \ge 2$.

\end{thm}

The constants $\gb_k$ are absolute, in particular,  independent of $C$. The first few values
of $\gb_k$ are:
$$
 \beta _2=0.373365,\ \beta _3=0.913728,\ \beta _4=1.52961, \beta_5=2.19252, \dots, \beta _{10}=5.8754, \dots
$$
Note that  the expected size\footnote{e.g. in the normal order sense of the  Erd\H{o}s-Kac theorem.} of $\gW(m)$ for a random integer $m$  is  $\log\log m$, 
and of course 
$$
\gW(x_{1,n}\cdot x_{2,n}\cdots x_{k,n})=\sum_j\gW(x_{j,n}),
$$ 
whence the expected size of 
this sum 
is 
$k\log\log C^n\sim k\log n$. Thus we may interpret \eqref{eq:gWliminf} as showing that, up to a multiplicative constant $k/\gb_k$, one never sees (asymptotically) a deficient number of prime factors.

To test the validity of this model empirically, it will be useful to understand how large $n$ should be to experimentally observe the behavior \eqref{eq:gWliminf}. Naively we may expect from this equation that the largest $\fn=n_{max}$ for which $x_{1,n}\cdots x_{k,n}$ is $R$-almost prime satisfies:
$$
{R\over \log \fn}\approx \gb_k,
$$
or 
\be\label{eq:naiveNmax}
\fn
 \approx \exp(R/\gb_k).
\ee
It turns out that  the probabilistic model sometimes makes a different prediction.

\begin{thm}\label{thm:main2}
Fix $k\ge 2$,
 $C>1$, and
for each $n\ge1$, draw a vector
$$
\bx_n=(x_{1,n},x_{2,n},\dots,x_{k,n})\in[1,C^n]^k
$$
uniformly. Let $\cX=(\bx_1,\bx_2,\dots)$ be a random variable consisting of a sequence of
 independent such draws, one for each $n$.
For any fixed $R\ge k$, consider the random variable
$$
\fn=\fn(R;\cX):=\max\{n:\gW(x_{1,n}\cdots x_{k,n})\le R\},
$$
with $\fn=0$ if there are no such $n$, and $\fn=\infty$ if the event occurs infinitely often. 
Then 
\begin{enumerate}
\item
with probability one,
\be\label{eq:main21}
\fn<\infty,
\ee
and moreover,
\item
for all   $m\ge k-1$, the $m$-th moment of $\fn$ diverges,
\be\label{eq:main22}
\bE[\fn^m]=\infty.
\ee
\end{enumerate}
\end{thm}

\begin{rmk} 
In the case  $k=1$ not covered in Theorem \ref{thm:main2}, one
has instead that with probability one, $\fn=+\infty$. 
\end{rmk}

\begin{rmk}\label{rmk:sporadic}
In many natural examples treated below, we have $k=2$, so taking $m=1$ means that the expected value of $\fn(R)$ is infinite  for all $R\ge2$. Thus we should not expect $\fn(R)$ to behave nicely like $\exp(R/\gb_2)$, as suggested naively by \eqref{eq:naiveNmax}. One may interpret this as saying that for $k=2$ there may exist extremely large ``sporadic'' solutions to $\gW(x_{1,n},\dots,x_{k,n})=R$. 
\end{rmk}

\begin{rmk}
The proofs of \thmsref{thm:main1} and \ref{thm:main2} apply and give the same result in the more general case of $\bx_n$ chosen from non-identically growing intervals, that is $
(x_{1,n},...,x_{k,n})\in[1,C_1^n]\times[1,C_2^n]\cdots\times[1,C_k^n]$,
for fixed constants $C_1,\dots,C_k>1$.
\end{rmk}

\subsection{The Toral Affine Sieve Conjecture}\

The  probabilistic model above,
motivates a heuristic prediction  concerning  the number of prime factors of certain sequences,
associated to toric orbits, 
the (rank one) ``Toral Affine Sieve Conjecture'' stated below. 
We will   derive as consequences of this conjecture other predictions
 in several settings of classical interest.

\begin{conj}[Toral Affine Sieve Conjecture]\label{conj:orb}
Let $\g\in\GL_2(\Q)$ be a hyperbolic matrix, that is, one having two distinct real eigenvalues; equivalently
$$
\tr(\g)^2-4\det(\g)>0.
$$ 
Let $\G=\<\g\>^+:=\{\g^n:n\ge0\}$ be the  semigroup generated by $\g$, and 
suppose that   $\bv_0\in\Q^2\setminus (0,0)$ is a nonzero vector
such that the orbit 
$\cO:=\G\cdot\bv_0\subset\Z^2$ is integral
and infinite. Then
\be\label{eq:orb}
\liminf_{(x,y)\in\cO}{\gW(xy)\over \log\log|xy|} \ \ge \ \gb_2 \approx0.373365 .
\ee
\end{conj}

Since  the Zariski closure of $\G$ in $GL(2, \C)$ is an algebraic torus, 
and since the orbit $\cO$ is assumed to be infinite, 
it is a one-dimensional torus, so it follows that
the Zariski closure of $\cO$ in $\C^2$  has
 $\dim(\Zcl(\cO))=1$ in \eqref{eq:Zcl}. 
We have taken the  test function $f(x,y) = xy$, whence $V(f)\cap\Zcl(\cO)$ is finite, having dimension $0$.
The points in $(x_n,y_n):=\g^n\bv_0\in\cO$ grow exponentially, that is, there are $C>c>1$ so that
$$
c^{n}\  < \  |x_ny_n| \ = \ |f(\g^n\bv_0)| \  <\ C^n
.
$$
In consequence,  the factor  $\log\log|xy|$ in  \eqref{eq:orb} can be replaced by $\log n$, that is,  \eqref{eq:orb} is equivalent to
$$
\liminf_{n\to\infty}{\gW(
x_n y_n
)\over \log n} \ \ge \ \gb_2 .
$$
 The conjecture is based on  applying the model of Theorem \ref{thm:main1} with $k=2$
 having two ``independent" factors $(x_{n}, y_{n})$ for $f(\g^n \bv_0)$.
In the ``generic" situation, we might have equality in these limits. However
there are cases of orbits whose limiting values may involve $\beta_k$ for larger $k$,
see the examples in \secref{sec:examples}.

\begin{rmk}
We did not need to assume in \conjref{conj:orb} any coprimality condition (e.g. $\gcd(\cO)=1$) on the orbit. Indeed, if all entries of $\bv=(x,y)\in\cO$ have a common factor, then
this factor, divided by $\log\log|xy|$, is irrelevant in the $\liminf$ in \eqref{eq:orb}.
\end{rmk}

\subsection{Consequences}\

The  basic \conjref{conj:orb} implies  other striking predictions, of which we present  two below; 
the first applies to integer points on affine quadrics, and the second applies to the
continued fraction convergents of quadratic surds.

\begin{thm}\label{thm:form}
Let 
$
Q(x,y)=Ax^2+Bxy+Cy^2
$ 
be an indefinite (that is, $D=B^2-4AC$ is positive), non-degenerate ($D$ is not a square) binary quadratic form over $\Z$. Fix a square-free $t\in\Z$ so that the set  $V(\Z)$ of $\Z$-points of the affine quadric $V=V_{Q,t}$
given by
$$
V\ : \ 
Q(x,y)
=t 
$$
is non-empty.
Then, 
assuming \conjref{conj:orb},
$$
\liminf_{{(x,y)\in V(\Z)}\atop{|xy| \to \infty}} 
{\gW(
xy
)\over \log\log|xy|} \    \ge \  \gb_2 .
$$
\end{thm}


\begin{thm}\label{thm:surd}
Let $\ga$ be a real quadratic irrational, and let $p_n/q_n$ denote the $n$-th  convergent
of its ordinary continued fraction expansion. 
Then,
assuming \conjref{conj:orb},
$$
\liminf_{n} {\gW(
p_nq_n
)\over \log n} \    \ge \  \gb_2. 
$$
\end{thm}

These two theorems will not be surprising to experts, but the (conditional) conclusions, particularly the appearance of the precise number $\gb_2\approx 0.373365$, are 
 unexpected.

\subsection{Organization}\

In  \secref{sec:examples}, we give a number of illustrative  examples and numerics which, one may argue, 
provide 
support for  the heuristic provided by the probabilistic model in the context of
 \conjref{conj:orb}.
We prove \thmref{thm:main1} in \secref{sec:pf1}, followed by \thmref{thm:main2} in \secref{sec:pf2}.
In the final \secref{sec:pf3}, we sketch proofs of \thmsref{thm:form} and \ref{thm:surd}.

\subsection{Notation}\

We use the following standard 
notation. We use the symbol $f\sim g$ to mean $f/g\to1$. The symbols $f\ll g$ and $f=O(g)$ 
are used 
interchangeably to mean the existence of an implied constant $C>0$ so that $f(x)\le C g(x)$ holds for all $x>C$; moreover $f\asymp g$ means $f\ll g\ll f$. 
Unless otherwise specified,  implied constants  
depend at most on $k$, which is treated as fixed.
The letter $\vep>0$ is an arbitrarily small constant, not necessarily the same at each occurrence. The Gamma function is denoted $\G(z)$ and a product $\prod_p$ denotes a product over primes. The floor function, $\lfloor\cdot\rfloor$, returns the largest integer not exceeding its argument.

\subsection*{Acknowledgements}\

The authors thank Jonathan Bober, Andrew Granville, Peter Sarnak, and Alireza Salehi Golsefidy for enlightening discussions, comments, and suggestions, and most of all, Danny Krashen and Sean Irvine for the highly non-trivial and time-consuming task of computing $\gW$ for Lucas, Fibonacci, and Mersenne numbers from cumbersome online databases of their factorizations.



\section{Examples and Numerics}\label{sec:examples}

It should be clear that
running decent numerics
to test \conjref{conj:orb}
is a daunting task. Indeed,  orbits increase exponentially in size, and hence become ever more difficult to factor. 
Thankfully, others have already exerted tremendous effort in tabulating prime factorizations 
for certain sequences of classical interest, in particular, the Fibonacci, Lucas, and Mersenne numbers. 
We  mine their factorization data to test our predictions for \conjref{conj:orb} and its consequences.
We have made the raw data and Mathematica file used to construct the figures available at: \url{http://sites.math.rutgers.edu/~alexk/files/AllOmegasData.nb}.

\subsection{Fibonacci and Lucas Numbers Factorization Statistics}\

\begin{figure}
\includegraphics[width=4in]{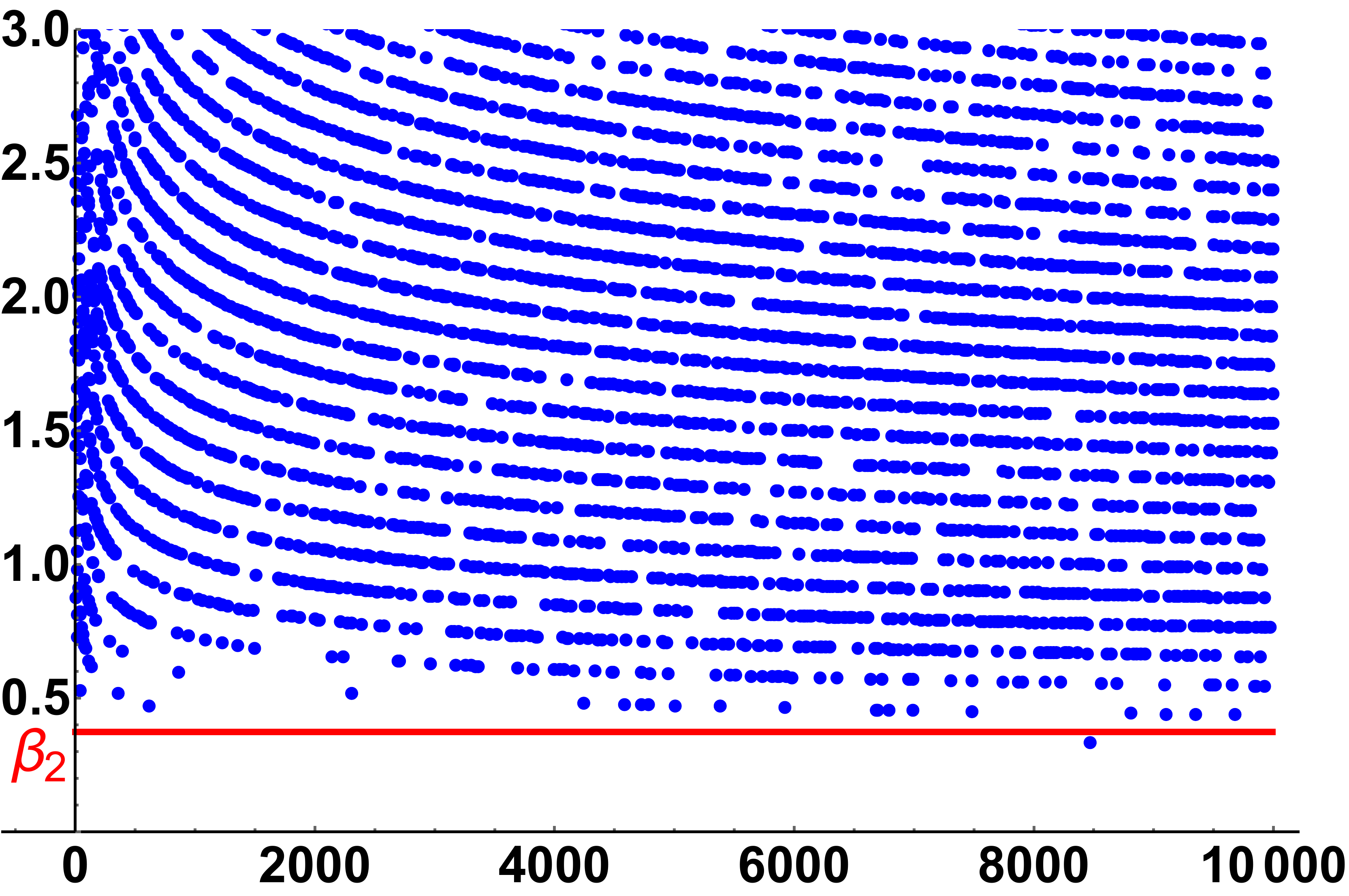}
\caption{A plot of $n<10\,000$ vs. $\gW(F_nL_n)/\log\log(F_nL_n)$. Also shown is the horizontal line $y=\gb_2\approx0.37.$}
\label{fig:1}
\end{figure}

Let $F_n$ and $L_n$ denote the $n$th Fibonacci and Lucas numbers, respectively. Recall that both sequences are defined by the same recursive relation, 
$F_{n+1}=F_{n}+F_{n-1}$ 
and 
$L_{n+1}=L_{n}+L_{n-1}$,
but differ in the initialization, namely, $F_1=F_2=1$, while $L_1=1$, $L_2=3$. They are related by 
\be\label{eq:F2nFnLn}
F_{2n}=F_nL_n.
\ee
Both sequences have been completely factored for $1\le n\le 1\,000$ and partially 
factored  for $n$ going  up to $10\,000$, see the website \cite{MersenneWeb}.

In the following calculations, when we encounter in the (incomplete) factorization data
a composite number having no known prime factors,  we treat that number as a product of exactly two primes  
(which may be an undercount in $\gW$).
We use this data to  study orbits giving several different combinations of Fibonacci numbers
and Lucas numbers.

\begin{Ex}\label{ex:FnLn}
\emph{One can easily verify that, if one takes 
$$
\g \ = \ \mattwo{1/2}{1/2}{5/2}{1/2}, \qquad
\G=\<\g\>^+,\qquad \bv_0=(1,1)^t, 
$$
then the orbit $\cO=\G\cdot\bv_0=\{(F_n,L_n):n\ge1\}$.
A plot of  $n$ versus
\be\label{eq:FLloglog}
{\gW(F_nL_n)
\over \log \log(F_nL_n)}
\ee
appears in \figref{fig:1}.  
This plot seems to give rather good evidence for equality in \eqref{eq:orb}. 
}
\end{Ex}

{\bf Remarks:}

$(i)$
The plot in \figref{fig:1} appears to be a union of curves, and a moment's thought reveals that these are roughly the level sets of $y=R/\log x$ for various integer values of $R$. Conjecture \ref{conj:orb}
predicts that the number of elements on each curve is finite, since each curve eventually dips below 
the line $y=\beta_2$.

$(ii)$
From \figref{fig:1}, one notices a single value of $n<10\,000$ for 
which 
\eqref{eq:FLloglog}
seems to dip below $\gb_2\approx0.37$. This occurs at $n=8\,467$, for which $L_n$ is prime and $F_n$ is composite, with each number spanning $1\,770$ decimal digits. 
Since we do not know any factors of $F_n$,  we follow our protocal, declaring that $\gW(F_nL_n)=3$. But the true value could perhaps be higher, in which case there may be no values of $n$ up to $10\,000$ dipping below \eqref{eq:FLloglog}. Since \conjref{conj:orb} only predicts a $\liminf$, there may in fact be infinitely many points in the plot dipping below $\gb_2$, as long as the amount by which they dip below decreases.

$(iii)$
The  data in  \figref{fig:1} also provide an instance of (the conditional) \thmref{thm:form}, since the pair
 $(F_n,L_n)$ are integer solutions to the Pellian binary quadratic form
\be\label{eq:Pell}
x^2-5 y^2=\pm4.
\ee


$(iv)$
While \figref{fig:1} may seem promising towards \conjref{conj:orb}, this computation is limited to
the humble scale $n=10\,000$, where $\log n\approx
\log\log(F_nL_n)\approx
10
$. 

With current  computing technology it would be difficultto go significantly 
farther.
\\

One may also object to using the Fibonacci and Lucas sequences to test \conjref{conj:orb}, as these are ``strong divisibility sequences''; i.e., $m\mid n\Longrightarrow a_m\mid a_n$. While it seems likely that this fact could affect some statistics of total number of primes seen in individual draws (see, e.g., \cite{BugeaudLucaMignotteSiksek2005}), it appears not to affect the $\liminf$ value in \eqref{eq:FLloglog}. Either way, any effect would only increase the limiting value, which \figref{fig:1} suggests is not the case.

\begin{figure}
\includegraphics[width=4in]{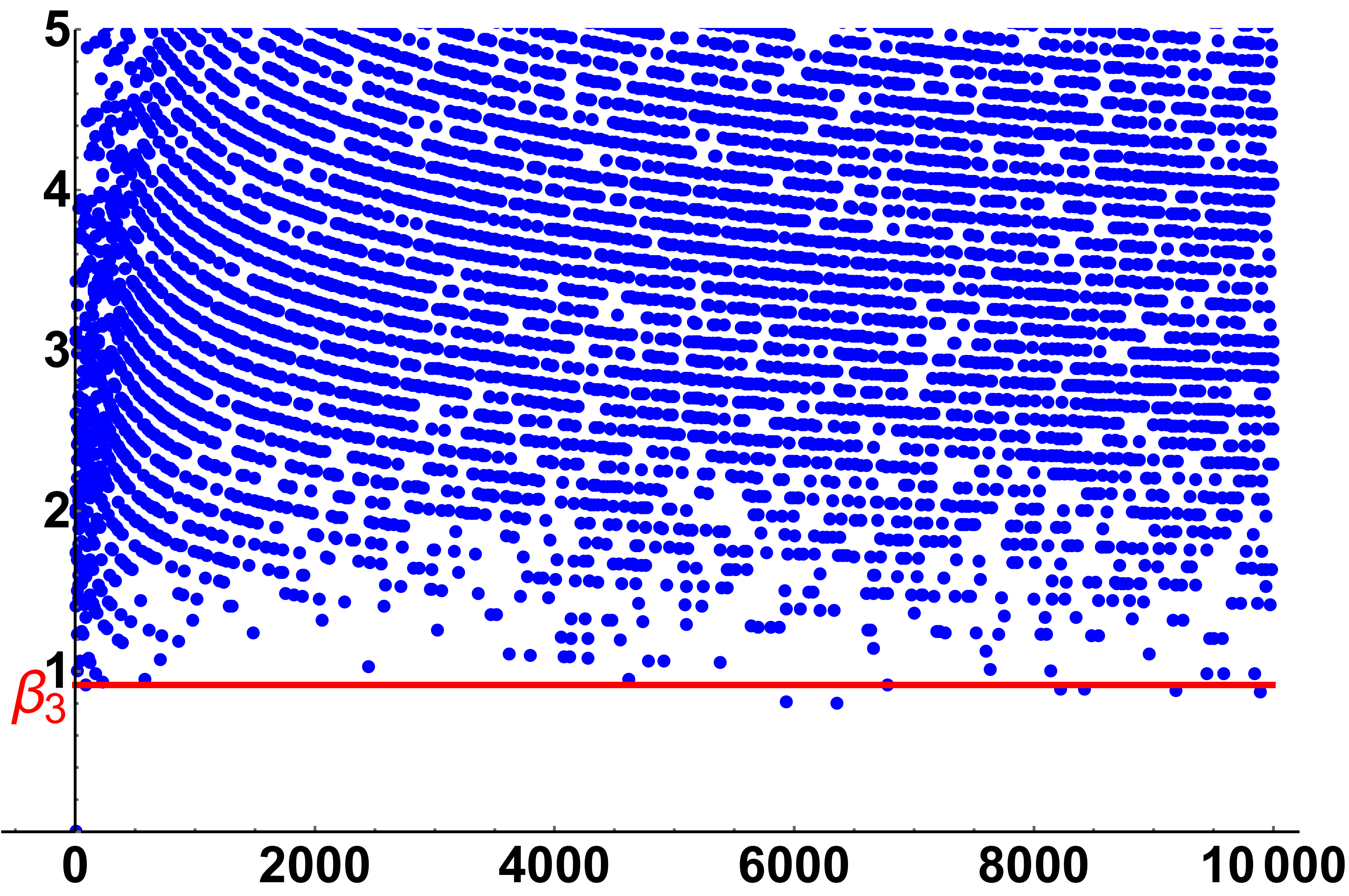}
\caption{A plot of $n<10\,000$ vs. $\gW(F_nF_{n+1})/\log\log(F_nF_{n+1})$. Also shown is the horizontal line $y=\gb_3\approx0.91.$}
\label{fig:2}
\end{figure}

\begin{figure}
\includegraphics[width=4in]{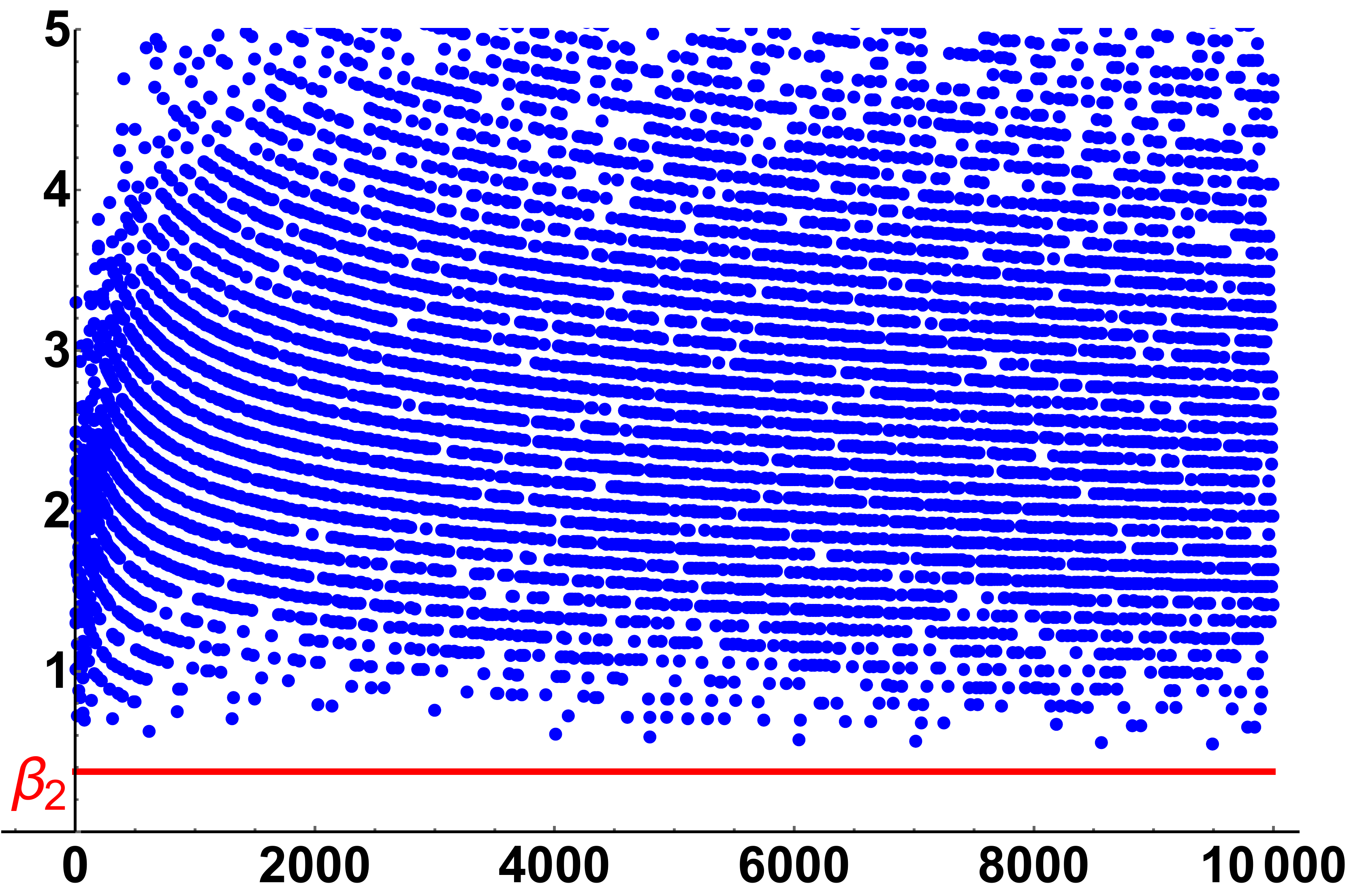}
\caption{A plot of $n<10\,000$ vs. $\gW(L_nL_{n+1})/\log\log(L_nL_{n+1})$. Also shown is the horizontal line $y=\gb_2.$}
\label{fig:3}
\end{figure}

\begin{Ex}
\emph{
Next we consider the simpler setting of 
consecutive Fibonacci numbers:
$$
\g \ = \ \mattwo 1110, \qquad
\G=\<\g\>^+,\qquad \bv_0=(1,0)^t, 
\qquad
\cO\ = \ \G\cdot\bv_0 \ = \ \{(F_{n+1},F_n)^t\}.
$$
Applying  \conjref{conj:orb}, one may surmise  that the correct liminf for $\gW(F_nF_{n+1})/\log\log(F_nF_{n+1})$ is $\gb_2\approx 0.37$. But a moment's inspection of \figref{fig:2} reveals that the truth seems to be closer to $\gb_3\approx0.91$. This is because one of the indices $n$ or $n+1$ is {\it even}, so that Fibonacci number splits according to \eqref{eq:F2nFnLn} into a Fibonacci times a Lucas. Thus this sequence $F_nF_{n+1}$ behaves like the product of {\it three} independent sequences, resulting in the predicted lim-inf of $\gb_3$, not $\gb_2$. 
}

\emph{
For this reason, \conjref{conj:orb} must be stated with an inequality in \eqref{eq:orb}; one cannot necessarily determine {\it a priori} from the data of $\cO$ whether there is a ``non-obvious'' factorization. Indeed, if we keep 
$\G$ as is but change $\bv_0$ to $\bv_0=(1,2)^t$, then the orbit $\cO=\{(L_{n+1},L_n)^t\}$ becomes consecutive Lucas numbers instead of Fibonaccis. These do not exhibit the extra factorization, so the liminf is restored (though now not very convincingly) to $\gb_2$, see \figref{fig:3}.
}
\end{Ex}

\begin{figure}
\includegraphics[width=4in]{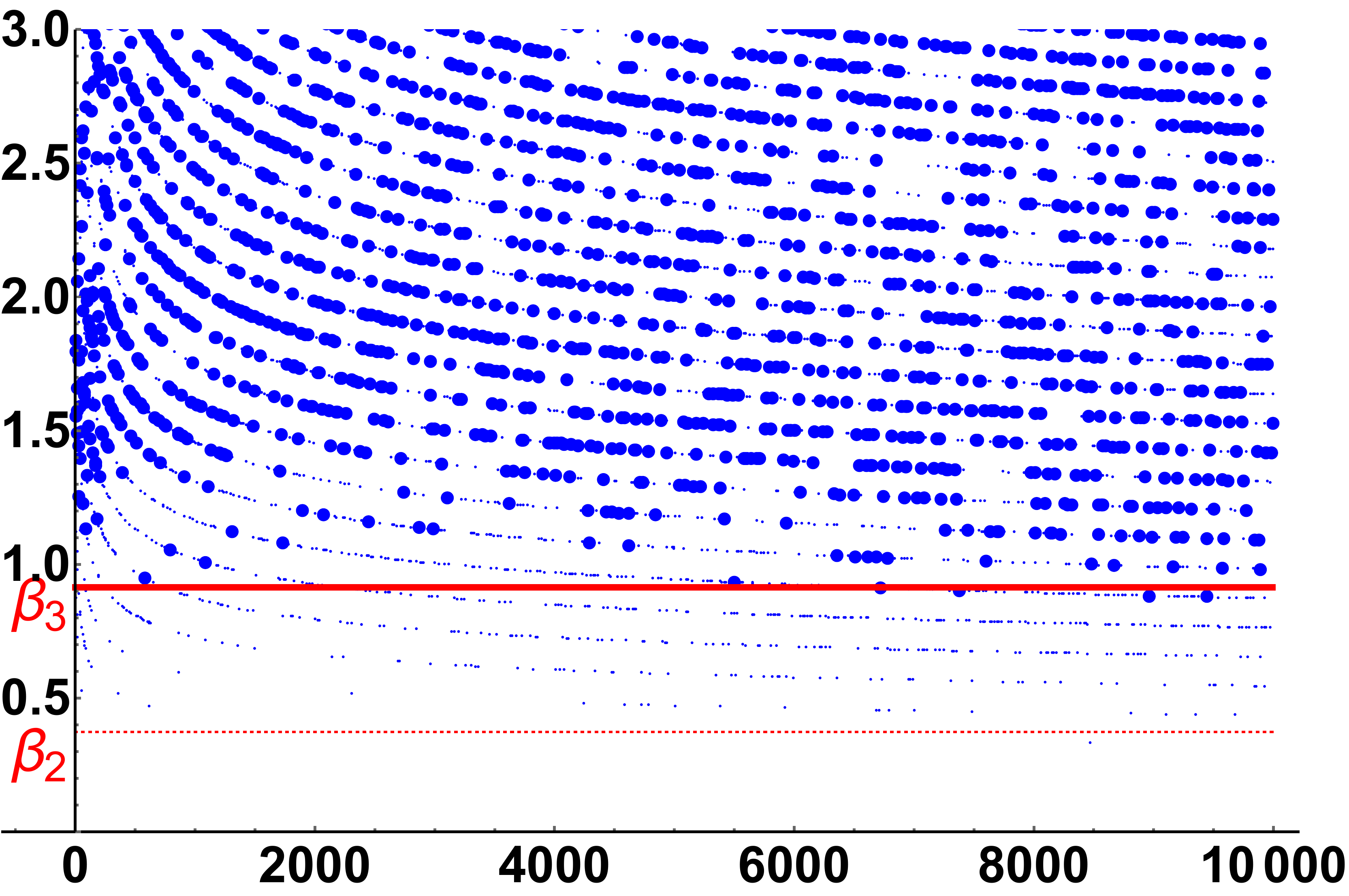}
\caption{A plot of $n<10\,000$ vs. $\gW(F_nL_{n})/\log\log(F_nL_{n})$, with the
even index values with large marks  and the odd index values with small marks. Also shown are the horizontal lines $y=\gb_2,\gb_3.$ Compare to \figref{fig:1}.}
\label{fig:FLeo}
\end{figure}

\begin{Ex}
\emph{
The previous example suggests the following refinement of \exref{ex:FnLn}. One can easily produce orbits which separately capture the even and odd index Fibonacci/Lucas pairs  $(F_{2n},L_{2n})$ and $(F_{2n+1},L_{2n+1})$. These of course appear simultaneously inside
the orbit of  \figref{fig:1}. Now  in \figref{fig:FLeo} we show what happens if the odd values are suppressed: the even values  exhibit an increased beta-value, again to $\gb_3$. 
}
\end{Ex}

\begin{figure}
\includegraphics[width=4in]{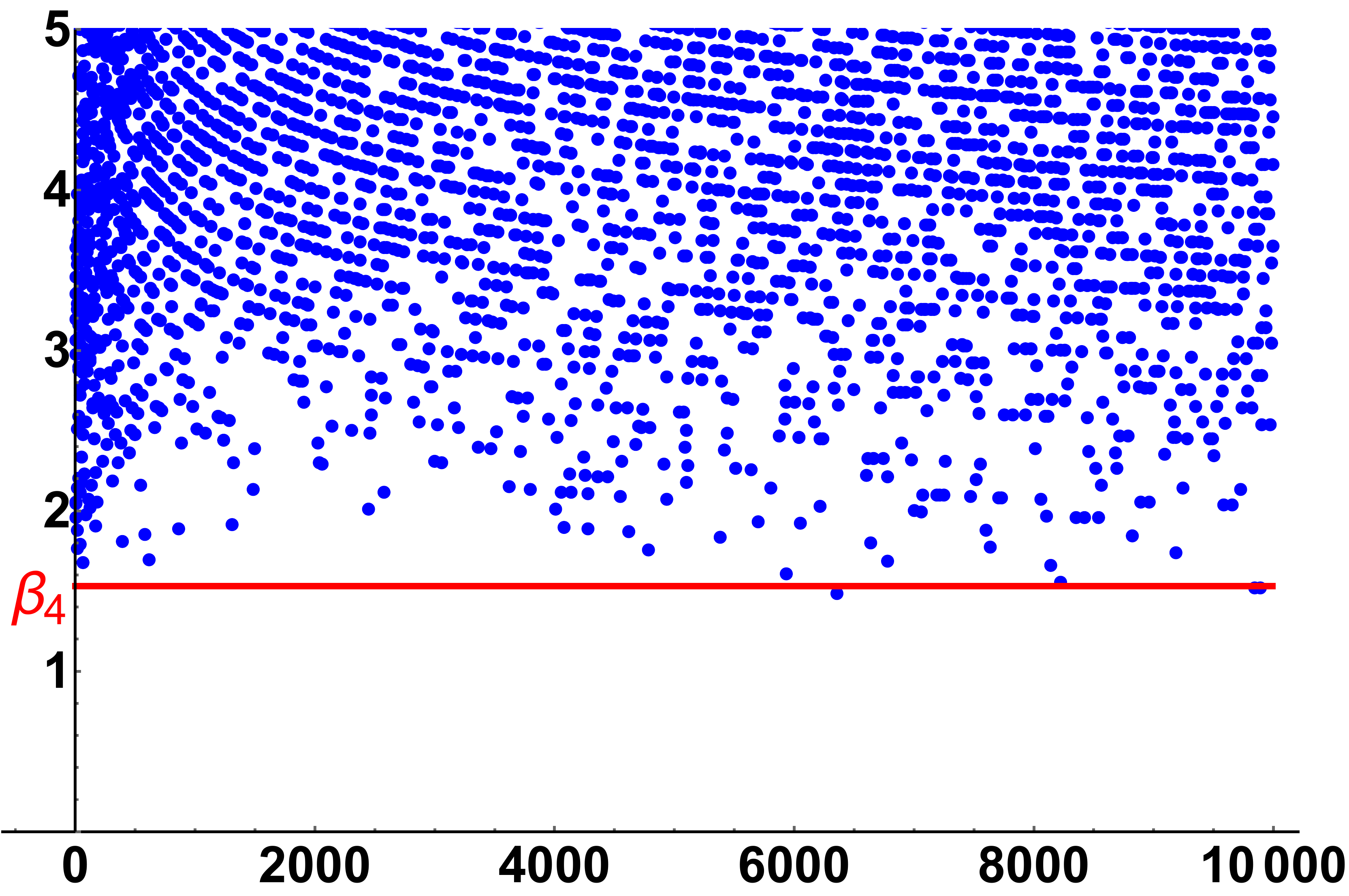}
\caption{A plot of $n<10\,000$ vs. $\gW(F_{2n}F_{2n+2})/\log\log(F_{2n}F_{2n+2})$. Also shown is the horizontal line $y=\gb_4.$}
\label{fig:F2nF2nm2}
\end{figure}

\begin{Ex}\label{Ex:trouble}
\emph{
We consider  pairs $(F_{2n},F_{2n+2})$ of consecutive even-indexed Fibonacci numbers.  This sequence was already discussed
 in the initial Bourgain-Gamburd-Sarnak paper on the Affine Sieve, see \cite[Section 2.1]{BourgainGamburdSarnak2010}.
It is obtained by taking 
$\g = \mattwo
3 1
{-1} 0 
$%
, which has powers
$$
\g^n = \mattwo
3 1
{-1}0 
^n  = \mattwo
{F_{2n+2}}  {F_{2n}}
{-F_{2n}} {- F_{2n-2}}
,
$$
and acting on  $\bv_0 = (1,0)^t$ to give the orbit $\cO=\{(F_{2n},F_{2n+2})^t\}$.
 Then 
 $$
 f( \g^n \bv_0) = F_{2n} F_{2n-2} = F_{n} L_{n} F_{n-1}L_{n-1},
$$
where we have again invoked the Fibonacci identity \eqref{eq:F2nFnLn}.
As a consequence we expect four ``independent" factors, so the liminf in \eqref{eq:orb} should be 
no smaller than $\gb_4 \approx 1.52961.$
See \figref{fig:F2nF2nm2}, which confirms the prediction. But on further inspection, it turns out that the lim-inf here should be $\gb_5$, not $\gb_4$! 
Indeed, one of the indices $n$ or $n-1$ is even, so one of the factors $F_n$ or $F_{n-1}$ in  $f(\g^n\bv_0)$ should always decompose further into a Fibonacci/Lucas pair. We do not fully understand why the numerics do not agree with this prediction, though it is plausible that the under-estimation of $\gW$ in inconclusive factorizations may at this point be making a significant contribution.
}
\end{Ex}

\subsection{Mersenne Number Factorization Statistics}\

For our last numerical example, we move to Mersenne numbers, $M_n:=2^n-1$, whose factorizations have also been extensively mined.

\begin{figure}
\includegraphics[width=4in]{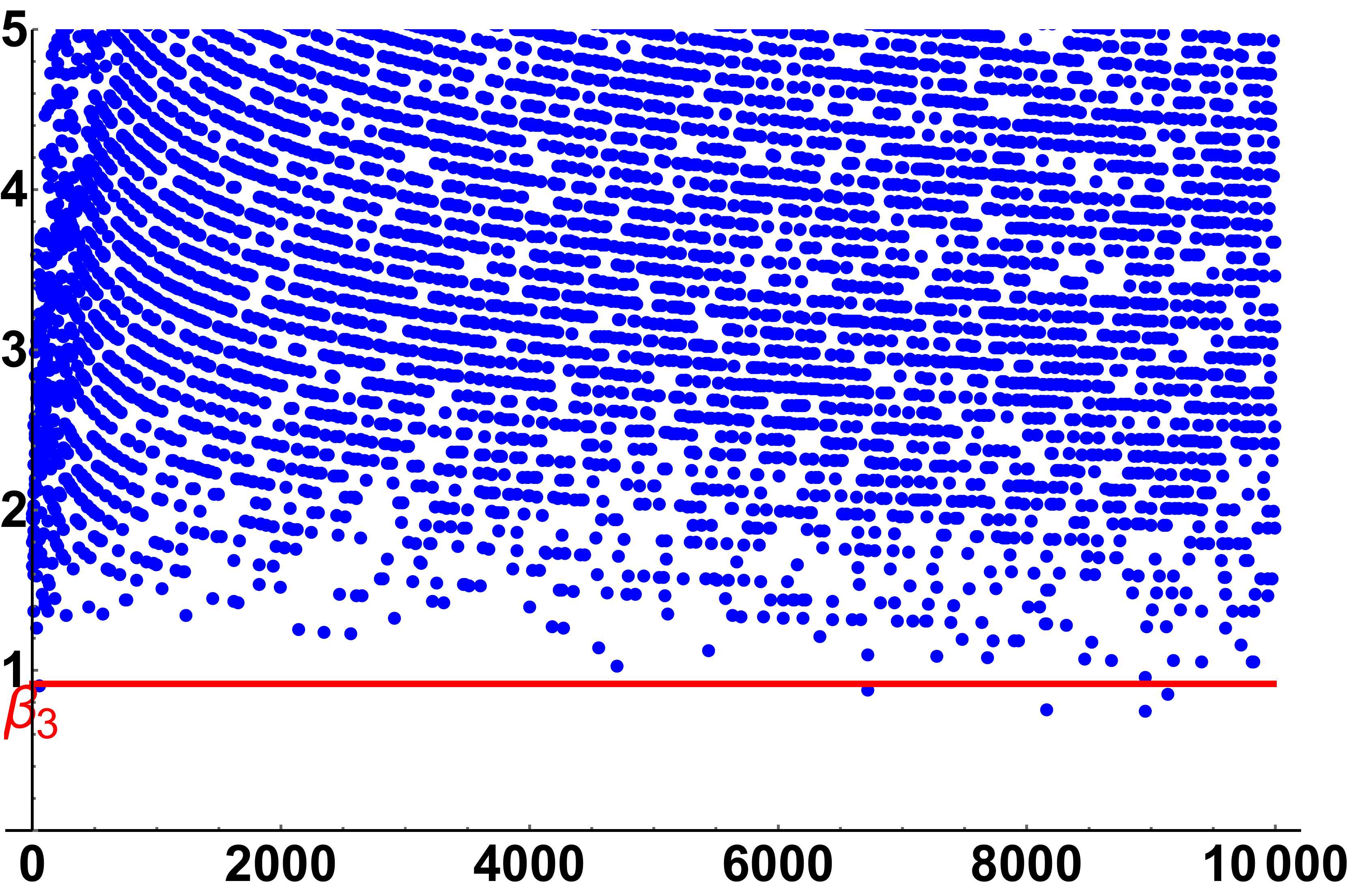}
\caption{A plot of $n<10\,000$ vs. $\gW(M_{n}M_{n+1})/\log n$. Also shown is the horizontal line $y=\gb_3.$}
\label{fig:6}
\end{figure}

\begin{Ex}
\emph{
To produce the orbit $\cO=\{(M_{n+1},M_n)\}$, consider as before $\G=\<\g\>^+$ and $\cO=\G\cdot\bv_0$, where:
$$
\g=\mattwo 3{-2}10, \quad \bv_0=(1,0)^t, \qquad \g^n\bv_0=(M_{n+1},M_n)^t.
$$
The first 500 values of $\gW(M_n)$ appear in OEIS (A046051), and the (sometimes partial) factorizations up to $10\,000$ were kindly provided to us by Sean Irvine using \url{factordb.com}. These were used to make \figref{fig:6}, showing that the liminf of $\gW(M_{n}M_{n+1})/\log \log (M_{n}M_{n+1})$ appears to be tending towards $\gb_3$. This is consistent with the fact that one of $n$ or $n+1$ is even, and for the even indices, 
Mersenne numbers $M_{2\ell}$ factor as $2^{2\ell}-1=(2^\ell-1)(2^\ell+1)$.
}
\end{Ex}

\subsection{Extreme Fibonacci and Lucas values with  a fixed number of prime factors}\


Let us now consider \thmref{thm:main2} and the (na\"ive) heuristic \eqref{eq:naiveNmax} in the case of the Fibonacci and Lucas sequences, for fixed $R=2$.

\begin{Ex}
\emph{
 Define the set
$$
\gS_{FF} \ := \ \{ n\ge2:\gW(F_nF_{n+2})=2\}
$$
to be the indices $n$ for which $F_n$ and $F_{n+2}$ are simultaneously prime.
Applying \eqref{eq:naiveNmax}
 with $R=k=2$ would suggest that
\be\label{eq:maxSig}
\max \gS_{FF} \ \overset{?}{\approx} \ \exp(2/\gb_2)\approx 212.
\ee
One can now examine the sequence \cite{OEISA001605}   of $n$  for which $F_n$
are prime, to find that 
\be\label{eq:ConsecF}
\{3, 5,11, 431, 569 \} \ = \ \gS_{FF} \cap [1,\, 1\, 000\,000].
\ee
Similarly,  consider the set
$$
\gS_{LL} \ := \ \{ n\ge2:\gW(L_nL_{n+2})=2\}
$$
of indices $n$ for which  $L_n$ and $L_{n+2}$ are simultaneously prime; presumably \eqref{eq:maxSig} should also hold for $\gS_{LL}$.
As before, one can examine the sequence \cite{OEISA001606}   of $n$  for which $L_n$
are prime, to find that 
\be\label{eq:ConsecL}
\{2, 5,11, 17  \} \ = \ \gS_{LL} \cap [1, \,1\,000\,000].
\ee
Both these results are compatible, at least to first order, with the naive heuristic \eqref{eq:maxSig}.
}
\end{Ex}

\begin{Ex}
\emph{
 Next define
$$
\gS_{FL} \ := \ \{ n\ge2:\gW(F_nL_n)=2\}
$$
to be the indices $n$ for which the Fibonacci and Lucas sequences are simultaneously prime.
As above, the naive heuristic   \eqref{eq:naiveNmax} predicts
$
\max\gS_{FL} \ \overset{?}{\approx} \ \exp(2/\gb_2)\approx 212.
$
Using  the sequences \cite{OEISA001605}  and  \cite{OEISA001606} of $n$ for which $F_n$ and $L_n$ are 
primes, respectively, however we find
\be\label{eq:gSis}
\{4, 5, 7, 11, 13, 17, 47, 148\,091 \} \ \overset{*}{=} \ \gS_{FL} \cap [1, \, 1\, 000\,000].
\ee
The ``$*$'' here is to note that for the largest index $\fn:=148\,091$, the corresponding $F_\fn$ and $L_\fn$ (each having around 30\,000 decimal digits) have not been certified prime.\footnote{The probable primality of $F_\fn$ was found by T. D. Noe while that of $L_\fn$ by de Water; see OEIS for further credits. 
Both numbers have passed numerous pseudoprimality tests.
Assuming GRH, one would need to run about $(30\,000)^4$ trials (that is, $(\log F_\fn)^2$ tests at a cost of $(\log F_\fn)^2$ each, ignoring epsilons) of the Miller primality test to certify these entries prime. Unconditionally, the exponent $4$ would be replaced by a $6$, see \cite{LenstraPomerance2011}. Or better yet, 
one could try the elliptic curve primality test, which is also unconditional and 
in practice runs faster, though
a worst-case execution time is currently unknown.}
The pair $(F_\fn,L_\fn)$, if indeed both entries are prime, would
have 
$$
{\gW(F_\fn,L_\fn)\over\log\log(F_\fn L_\fn)} \ \overset{?}{\approx} \ \frac2{\log \fn}  \ \approx\ 0.167988,
$$
so if we extended \figref{fig:1} to $n<150\,000$,  we would see a huge dip below $\gb_2$ at $\fn$.
In light of \eqref{eq:maxSig}, this certainly constitutes a massively ``sporadic'' solution to \eqref{eq:Pell}.
However  but the existence of such a solution  is {\it not} shocking, as it is predicted
to sometimes occur by the probabilistic model of \thmref{thm:main2}
(see \rmkref{rmk:sporadic}).
It seems likely to us (though again, this may be na\"ive) that the left side of \eqref{eq:gSis} is actually an equality
to $\gS_{FL}$.\footnote{Note that in some very special cases, one 
can
sometimes
 completely determine sets like $\gS_{FL}$. Indeed, see \cite{BoberLagariasSchmuland2009}, where all solutions to $x^2-3y^2=1$ with $\gW(xy)\le3$ are effectively listed.}
}
\end{Ex}


\section{Proof of \thmref{thm:main1}}\label{sec:pf1}

\subsection{Analysis of $\gb_k$}\

Fix  an integer $k\ge 1$  let $\gb_k$ solve \eqref{eq:gbkDef}. We first analyze this equation.

\begin{lem}\label{lem:fk-bound}
For real $k \ge 1$ the function 
$$
f_k(t):=t(1-\log t+\log k)-(k-1)
$$
is increasing on $0<t<k$. It has a unique root $t=\gb_k\in(0,k-1].$
\end{lem}
\pf
The derivative of $f$ is $f'_k(t)=-\log t+\log k,$ which is clearly positive on $(0,k)$. 
For $k=1$ it has by inspection a root at $\beta_0=0 = k-1$.
For $k >1$, near the origin,  
 $$
 \lim_{t\to0^+}f_k(t) \ = \ -(k-1) \ < \ 0,
 $$ 
 and at $t=k-1$, we have 
 $$
 f_k(k-1) \ = \ (k-1)\log\left(\tfrac k{k-1}\right) \ > \ 0.
 $$
Hence $f_k(t)$ has a unique root in this interval.
\epf
\begin{rmk}
One can 
solve for $\gb_k$ explicitly in terms of the inverse function $g(z)$  to $z\mapsto ze^z$ on the positive
real axis.  Namely, one finds
$$
\gb_k={1-k\over g\left({1-k\over ek}\right)},
$$
where $e=2.718...$.
We will not need this fact, nor the fact 
that $\gb_k=k-1-O(1/k)$ for $k$ large, which can be shown in a variety of ways.
\end{rmk}
\subsection{Analysis of the behavior of $\gW$}\

We next  record a uniform asymptotic formula for
$$
\cN_r(T)\ := \ \#\{x<T:\gW(x)=r\},
$$
that is,  the number of positive integers up to $T$ having exactly $r$ prime factors, counted with multiplicity.
For fixed $r$, the formula
\be\label{eq:naiveNrT}
\cN_r(T)\ \sim \ {T\over \log T}{(\log\log T)^{r-1}\over (r-1)!}
,\qquad\qquad(T\to\infty)
\ee
is well-known, but we shall require an estimate when $r$ is an increasing function of  $T$. 
Such an estimate can be obtained based on   a method of Selberg \cite{Sel54}. 
A treatment is given in  Tenenbaum \cite[Chap. II.6, Theorem 5]{Tenenbaum1995}, as stated below.

The result is given in terms of  the function
$$
\nu(z)\ := \, \frac{1}{\G(z+1)}\prod_p 
\left(
 \left(1-\frac zp\right)^{-1}\left(1-\frac1p\right)^z 
 \right).
$$
This infinite product converges on $\Re(z)>0$, giving in this region a non-vanishing meromorphic function with simple poles at $z=p$ for all primes $p$. Note also that $\lim_{z\to0^+}\nu(z)=1$; hence for real   $z\in [0,3/2]$, say, 
$\nu(z)$ is  bounded above and below by positive constants.

\begin{prop}[{\cite[eqn. (20), p. 205]{Tenenbaum1995}}]
For $T\ge3$, we have uniformly in 
$$
1\le r \le \frac32\log\log T
$$ 
that
\be\label{eq:Tenen}
\cN_r(T)\ = \ 
{T\over \log T}{(\log\log T)^{r-1}\over (r-1)!}\left(\nu\left({r-1\over \log\log T}\right)+O\left({r\over (\log\log T)^2}\right)\right)
,
\ee
with an absolute implied constant.
\end{prop}


This asymptotic continues to hold up to $r<(2-\gep)\log\log T$, but not beyond this point, as $\nu$ has a pole at $z=2$. A different asymptotic formula takes over at $r>(2+\gep)\log\log T$, see \cite{Nicolas1984}, but it will not be needed for our purposes.

For our application we derive from  \eqref{eq:Tenen}
a simplified estimate.

\begin{lem}
Let $r=\g\log\log T$ with $\frac1{\log\log T}\le\g<\frac32$. Then as $T\to\infty$,
\be
\label{eq:NcT} 
\bP[\gW(x)=r]  \ := \ \frac{\cN_r(T)}T
\ \asymp \
(\log T)^{\g-\g\log\g-1+o(1)}
,
\ee
with absolute implied constants. 
\end{lem}
\pf
First recall that, on $[0,3/2]$, the function $\nu(\cdot)$ is bounded above and below by positive constants. Then inserting the Stirling's formula estimate,
$$
(r-1)! \ \asymp \ r^{r-1/2}e^{-r},\qquad\qquad(1 \le r < \infty)
$$ 
into \eqref{eq:Tenen} yields
\beann
\bP[\gW(x)=r] &\asymp &
{1\over \log T}\left({\log\log T\over r}\right)^{r-1}r^{-\foh} e^r
\ = \
{1\over \log T}
\left({\g}\right)^{-\g\log\log T-1}
(\g\log\log T)^{-\foh} 
(\log T)^\g
\\
& = &
\g^{-\frac32}
(\log\log T)^{-\foh}
(\log T)^{-1-\g\log\g+\g}
,
\eeann
from which the estimate \eqref{eq:NcT} follows, since $\g\ge1/\log\log T$.
\epf

\subsection{Estimate for a single draw}\

To prove \thmref{thm:main1}, we first obtain 
upper and lower bounds on the probability density function for a single draw. 

\begin{thm}\label{thm:pf1}
Let $k \ge 1$ be fixed.
For any integer $T\ge2$, draw a vector 
$$
(x_1,x_2,\dots,x_k)\in [1,T]^k
$$
uniformly. For any small $\vep>0$, there is a $\gd=\gd(\vep)>0$ so that  for all $T>T_0(\vep)$,
\be\label{eq:pf12}
\bP[\gW(x_1x_2\cdots x_k)\le (\gb_k+\vep)\log\log T] \ \gg_\vep \ \frac1{(\log T)^{1-\gd}},
\ee
and,  for $k \ge 2$, 
\be\label{eq:pf11}
\bP[\gW(x_1x_2\cdots x_k)\le (\gb_k-\vep)\log\log T] \ \ll_\vep \ \frac1{(\log T)^{1+\gd}},
\ee
  there is a $\gd=\gd(\vep)>0$ so that  for all $T>T_0(\vep)$,
\end{thm}

\subsubsection{Proof of the lower bound \eqref{eq:pf12}}\
 
Suppose $k \ge 1$ and write $k\g=\gb_k+\vep$, so that $0<\g<1$, and 
let
$$
r \ := \ \left\lfloor
{\g}\log\log T
\right\rfloor.
$$
Then
$$
\bP[\gW(x_1\dots x_k)\le k\g\log\log T]
\ \ge\
\prod_{j=1}^k
\bP[\gW(x_j)=r].
$$
Inserting \eqref{eq:NcT}  gives
\beann
\bP[\gW(x_1\dots x_k)\le k\g\log\log T]
& \gg&
\left[
(\log T)^{\g-\g\log\g-1+o(1)}
\right]^k
.
\eeann
Write $\ga=k\g$; then as $T \to \infty$ the exponent of $\log T$
approaches the limiting value
$$
\g k-\g k \log\g-k
\ =\
\ga
-\ga \log\ga
+\ga \log k
-k  = f_k(\alpha) - 1
.
$$
By  Lemma \ref{lem:fk-bound}, since
$\ga=k\g=\gb_k+\vep>\gb_k$, and $f_k(\beta_k)=0$, we conclude
that as $T \to \infty$
the limiting exponent exceeds $-1$ by the positive  amount $f_{k}(\ga)  >0$.
Therefore we can pick $\delta(\epsilon) >0$ and $T _0(\epsilon)$ 
depending on $\vep$ (and $k$, which is fixed) so that \eqref{eq:pf12} holds.

\subsubsection{Proof of the upper bound \eqref{eq:pf11}}\

The upper  bound estimate \eqref{eq:pf11}   is more subtle and requires $k \ge 2$.
Again take a  fixed $\epsilon>0$  and  define $\gamma$ by $k\g=\gb_k-\vep$ 
taking $\epsilon$ small enough  that  $0<\g<1$, which is possible since $\gb_k >0$.
Since
$$
\gW(x_1\cdots x_k) \  = \ \gW(x_1)+\cdots+\gW(x_k),
$$
we have that
$$
\bP[\gW(x_1\dots x_k)\le k\g\log\log T]
\ = \
\sum_{r_1+\cdots+r_k\le  k \g\log\log T}
\bP[\gW(x_1)=r_1,\cdots \gW(x_k)= r_k ]
$$
We upper bound the total number of summands trivially by 
$$
\sum_{r_1+\cdots+r_k\le  k \g\log\log T}1 
\ \ll
\
(\log\log T)^k
\ = \ 
(\log T)^{o(1)}. 
$$
It  remains to upper bound the contribution of an individual summand 
$$
\max_{ r_1+\cdots +r_k \le  k \g\log\log T}\
\bP[\gW(x_1)=r_1,\cdots \gW(x_k)= r_k ]
.
$$
Write each $r_j$ as 
$$
r_j \  = \ \g_j\log\log T,
$$
so that
\be\label{eq:gjSum}
\g_1+\cdots +\g_k\  \le\ k\g\ <\ \gb_k \ <\ k-1.
\ee
On average these $\g_j$'s are less than one, but individually they could in principle be large, and we can apply \eqref{eq:NcT} only when  $\g_j<3/2$.
Let $\ell\subset\{1,\dots,k\}$ denote the indices $j$ for which $\g_j<3/2$ is ``low,'' and let $h:=\{1,\dots,k\}\setminus\ell$ be the ``high'' indices. Abusing notation, we use the same symbol for their cardinalities, e.g.,
 $$
 \ell+h=k.
 $$ 
 We have that 
$$
k\g \ \ge \ \sum_{j\in h}\g_j\ \ge \ \tfrac{3}{2}h,
$$
so
$$
\ell \ \ge\ k(1-\tfrac23\g) \ > \ \tfrac{1}{3}k,
$$
and
\be\label{eq:gjLow}
\sum_{j\in \ell}\g_j \ = \ \sum_{j}\g_j - \sum_{j\in h}\g_j \ \le \ k\g - \tfrac32h.
\ee
For $j\in h$, we estimate $\bP[\gW(x_j)=r_j]\le 1$ trivially. This gives a bound
\beann
\bP[\gW(x_1)=r_1,\cdots \gW(x_k)= r_k ]
& \le &
\prod_{j\in\ell}
\bP[\gW(x_j)=r_j]
\\
& \ll &
(\log T)^{o(1)}
\prod_{j\in\ell}
(\log T)^{\g_j-\g_j\log\g_j-1}
,
\eeann
using \eqref{eq:NcT}. The exponent in this
expression,  subject to \eqref{eq:gjLow}, is maximized if, for all $j\in\ell$,
we set all values equal $\g_j= \eta$, 
in which case,
\be\label{eq:bPup}
\bP[\gW(x_1)=r_1,\cdots \gW(x_k)= r_k ]
\ \ll \
(\log T)^{\ell(\eta-\eta\log\eta-1)+o(1)}
.
\ee
Now we have  
$$
\eta \ = \ \eta(\g,k,\ell) \ : = \
{k\g\over \ell}-{3h\over 2\ell}
\ = \
{3\over 2}
-{k\over \ell}
\left(
{3\over 2}
-{\g}
\right)
.
$$
We bound the exponent \eqref{eq:bPup}, varying $\ell$.
Viewing $\ell$ as a continuous variable, we 
$$
 \eta' \ := \ 
{\dd\eta\over\dd\ell} \ = \
{k\over \ell^2}
\left(
{3\over 2}
-{\g}
\right)
\ = \
\frac1\ell\left(\frac32-\eta\right).
$$
The derivative of the exponent of $\log T$ is in the $\ell$-variable is then
\beann
{\dd\over\dd\ell}[\ell(\eta-\eta\log\eta-1)] 
& = &
\eta-\eta\log\eta-1
-\ell
\eta'\log \eta
\\
& = &
\eta
-
\frac32\log \eta
-1
,
\eeann
which by inspection is a positive function of $\eta\in(0,1)$. 
It follows that  the exponent is maximized at the largest allowable  value of $\ell$, namely
the integer $\ell=k$, so $h=0$. 
For this value of $\ell$, we have $\eta=\g$, whence as $T \to \infty$ 
the exponent of $\log T$ in \eqref{eq:bPup} 
approaches  the limiting value
$$
k(\g-\g\log\g-1)
=
\ga-\ga\log\ga
+\ga\log k-k = f_{k}(\ga) -1.
$$
where we have again set $\ga=k\g=\gb_k-\vep$.
Again using  \lemref{lem:fk-bound}
this  limiting exponent 
 is less than $-1$  since  $\ga < \gb_k$ gives  $f_k(\ga) < 0$. 
Thus we can choose $\delta(\epsilon)$ and a $T_0(\epsilon)$ so that
\eqref{eq:pf11} holds.
 This completes the 
proof of \thmref{thm:pf1}.
\\

\subsection{Proof of \thmref{thm:main1}}\

It is now a simple matter to deduce \thmref{thm:main1} from \thmref{thm:pf1}.
Instead of a single  draw, here we have a sequence of independent draws, one for each $n=1,2,\dots$, and with $T=C^n$. By \eqref{eq:pf11},
$$
\bP\left[{\gW(x_{1,n}x_{2,n}\cdots x_{k,n})\over\log n}\le (\gb_k-\vep)(1+\log\log C/\log n)\right] \ \ll_\vep \ \frac1{n^{1+\gd}},
$$
and $\sum_{n\ge1}1/n^{1+\gd}<\infty$. Thus by the Borel-Cantelli Lemma, 
the probability of these events occurring infinitely often is zero; that is, with probability one, we have
$$
\liminf_n{\gW(x_{1,n}x_{2,n}\cdots x_{k,n})\over\log n}\ge \gb_k-\vep.
$$

Similarly, the independent events 
$$
\left[{\gW(x_{1,n}x_{2,n}\cdots x_{k,n})\over\log n}\le (\gb_k+\vep)(1+\log\log C/\log n)\right]
$$
occur with probability at least $1/n^{1-\gd}$, the sum of which diverges. By the second Borel-Cantelli Lemma,  infinitely many occur with probability one, so
$$
\liminf_n{\gW(x_{1,n}x_{2,n}\cdots x_{k,n})\over\log n}\le \gb_k+\vep.
$$
This proves \thmref{thm:main1}.
\\


\section{Proof of \thmref{thm:main2}}\label{sec:pf2}

Let $k\ge 1$, $C>1$, and $R\ge1$ be fixed throughout this section (unlike the previous section, where $R$ was growing).
In particular, the estimate \eqref{eq:naiveNrT} is perfectly valid here and will be used regularly. In this section, we allow implied constants to depend on $k, C$ and $R$, since they are fixed.

For each $n\ge1$, we choose uniformly a vector $\bx_n=(x_{1,n},\dots,x_{k,n})\in[1,C^n]^k$, and let 
$$
\fn \ = \ \fn(R) \ = \ \max\{n\ge1:\gW(x_{1,n}\cdots x_{k,n})\le R\},
$$
with $\fn=0$ if this set is empty and $\fn=\infty$ if it is unbounded.

First note that \eqref{eq:main21} follows immediately from \thmref{thm:main1}. Indeed, if $\fn(R)=\infty$, then $\gW(x_{1,n}\cdots x_{k,n})=R$ occurs for infinitely many $n$'s. But then 
$$
\liminf_{n\ge1}{\gW(x_{1,n}\cdots x_{k,n})\over \log n}\ = \ 0,
$$
contradicting \eqref{eq:gWliminf}. Hence this event has probability zero.

To prepare for the proof of \eqref{eq:main22}, we record the following computations.
Recall that  implied constants in this section may depend on $k$, $C$, and $R$.

\begin{lem}
Let $k \ge 1$ and $R \ge 1$ be fixed. Then for $t \ge 1$, 
\be\label{eq:bPupper}
\bP[\gW(x_{1,t}\cdots x_{k,t})\le R] \ \ll \
{(\log t)^{k(R-1)}\over t^k}.
\ee
Assuming further that $R\ge k$, we have that
\be\label{eq:bPlower}
\bP[\gW(x_{1,t}\cdots x_{k,t})\le R] \ \gg \
{(\log t)^{R-k}\over t^k}.
\ee
\end{lem}
\pf
The event $\gW(x_{1,t}\cdots x_{k,t})\le R$ is contained inside the intersection of the events $\gW(x_{j,t})\le R$, for all $j=1,2,\dots,k$. Thus using \eqref{eq:naiveNrT} gives
\beann
\bP[\gW(x_{1,t}\cdots x_{k,t})\le R] & 
\le &
\prod_{j=1}^k
\bP[\gW(x_{j,t})\le R] 
\
\ll
\
\left[
{1\over \log C^t}
{(\log \log C^t)^{R-1}
\over (R-1)!}
\right]^k
,
\eeann
from which \eqref{eq:bPupper} follows immediately.

Now assume that $R/k\ge 1$. Then 
the event $\gW(x_{1,t}\cdots x_{k,t})\le R$  contains the intersection over all $j=1,2,\dots,k$ of the non-empty events $\gW(x_{j,t})\le R/k$. So
\beann
\bP[\gW(x_{1,t}\cdots x_{k,t})\le R] & 
\ge &
\prod_{j=1}^k
\bP[\gW(x_{j,t})\le \tfrac Rk] 
\
\gg
\
\left[
{1\over \log C^t}
{(\log \log C^t)^{\frac Rk-1}
\over (\frac Rk-1)!}
\right]^k
,
\eeann
which implies \eqref{eq:bPlower}.
\epf

\begin{lem}\label{cor:bnd}
If $R\ge k\ge1$ are fixed, then for all sufficiently large $t$, 
$$
\bP[\fn(R)=t]
\ \gg\
{(\log t)^{R-k}\over t^k}
.
$$
\end{lem}

\pf
Consider the event $\fn(R)=t$. This occurs if and only if $\gW(x_{1,t}\cdots x_{k,t})\le R$ and, for all larger integers
 $s>t$, we have that $\gW(x_{1,s}\cdots x_{k,s})> R$. That is,
\beann
\bP[\fn(R)=t]
&=&
\bP[\gW(x_{1,t}\cdots x_{k,t})\le R]
\cdot
\prod_{s>t}
\bigg(
1-
\bP[\gW(x_{1,s}\cdots x_{k,s})\le R]
\bigg)
\\
&\gg&
{(\log t)^{R-k}\over t^k}
\cdot
\prod_{s>t}
\bigg(
1-
K{(\log s)^{k(R-1)}\over s^k}
\bigg)
,
\eeann
where we used \eqref{eq:bPlower} and \eqref{eq:bPupper}. (Here $K>0$ is a constant depending at most on $k$, $C$, and $R$.)
Since $s\ge2$, the infinite product converges absolutely. It bounds the result below by a uniform positive
constant for all sufficiently large $t$ that  avoid possible nonpositive  terms 
for small $s$ 
in the infinite product.
\epf

\pf[Proof of \thmref{thm:main2}]
Assume that $R\ge k\ge1$ and let $m\ge k-1$. Consider the $m$-th moment of $\fn$, namely,
\beann
\bE[\fn^m] & = &
\sum_{t\ge0} t^m\, \bP[\fn(R)=t]
\ \gg \
\sum_{t\ge0} t^m 
{(\log t)^{R-k}\over t^k}
,
\eeann
where we used \lemref{cor:bnd}.
Since $m-k\ge-1$, this sum diverges.

Note the case  $R=k=1$ gives divergence of the $m=0$-th moment; that is, if $k=1$ then  $\fn=\infty$
with probability $1$.)
\epf


\section{Proofs of \thmsref{thm:form} and \ref{thm:surd}}\label{sec:pf3}

Assume \conjref{conj:orb} in this section.
\pf[Proof of \thmref{thm:form}] 
Let $V:Q=t$ have $V(\Z)\neq\emptyset$.
As is well-known and in this case essentially goes back to Gauss, $V(\Z)$  decomposes into a finite number of $\G$-orbits, 
$$
V(\Z) \ = \ \bigsqcup_{j=1}^m \G\cdot \bv_j,
$$
where $\G=O_Q(\Z)$ is the orthogonal group fixing $Q$ (see, e.g., \cite{Cassels1978} or \cite[\S2]{Kontorovich2016}). Since $Q$ 
is indefinite,
 the Zariski closure of $\G$ is a torus, $$
\bbG=\Zcl(\G)=O(1,1).
$$ 
Thus,
 up to finite index, $\G=\<\g\>$ for some hyperbolic matrix $\g$. By \conjref{conj:orb} each orbit $\cO_j=\G\cdot \bv_j$ has 
$$
\liminf_{(x,y)\in\cO_j}{\gW(xy)\over \log\log|xy|} \ \ge\ \gb_2,
$$
and hence the same holds for all of $V(\Z)$.
\epf

\pf[Proof of \thmref{thm:surd}]
Let $\ga$ be a quadratic surd having
ordinary  continued fraction expansion $\ga=[a_0, a_1, a_2,...]$
with partial quotients $p_n/q_n$, given in matrix form by
$$
\mattwo 0 1 1 0\mattwo 0 1 1 {a_{0}}  \mattwo 0 1 1 {a_{1}}\cdots  \mattwo 0 1 1 {a_{n}} \vectwo01 = \vectwo {p_{n}} {q_{n}}.
$$
Now $\ga$ has an eventually periodic continued fraction expansion 
$$
\ga \ = \ [a_0;a_1,\dots a_k, \overline{a_{k+1},\dots,a_{k+\ell}}].
$$
After the first few terms, the sequence $(p_n,q_n)^t$ decomposes into finitely many $\G$-orbits, where 
$$
\G=\<\g\>,    \qquad \g=M \mattwo 0 1 1 {a_{k+1}}\cdots  \mattwo 0 1 1 {a_{k+\ell}}M^{-1}, 
$$
with 
$$
M= \mattwo 0 1 1 0\mattwo 0 1 1 {a_{0}}  \mattwo 0 1 1 {a_{1}}\cdots  \mattwo 0 1 1 {a_{k}},
$$
for the orbits given by 
$$                                      
 \bv_j := M \mattwo 0 1 1 {a_{k+1}}\cdots  \mattwo 0 1 1 {a_{k+j}} \vectwo 01, 
 \quad \quad 0 \le j \le \ell-1.
$$
We may apply  \conjref{conj:orb} to each orbit, since they are infinite,  and using the asymptotic  $\log\log p_n q_n \sim \log n$ establishes the result.
\epf

\end{document}